\documentclass[12pt,reqno,a4wide]{amsart}

\usepackage{tikz} 
\usepackage{calrsfs}
\usepackage[charter]{mathdesign}

\allowdisplaybreaks

\begin{document}
	
\title[Lattice paths, Elenas, and bijections]
{Height restricted lattice paths, Elenas, and bijections}
\author[Helmut Prodinger]{Helmut Prodinger}
\address{Department of Mathematics, University of Stellenbosch 7602,
	Stellenbosch, South Africa}
\email{hproding@sun.ac.za}
\thanks{The author was supported by an incentive grant of the National Research Foundation of South Africa.}


\begin{abstract}A  bijection is constructed between two sets of height restricted lattice paths by means of translating them in two tree classe, namely plane trees and Elena trees. An old bijection between them can be used now for that actual problem.
\end{abstract}

\date{}

\maketitle
\section{Introduction}

We consider lattice path, consisting of up-steps and down-steps of one unit each, starting at the origin. In other words,
$(0,s_0), (1,s_1),\dots ,(n,s_n)$, with $s_0=0$, and $|s_i-s_{i+1}|=1$. Furthermore, we assume that
$0\le s_i\le 3$. Let $\mathcal{A}_{n,i}$ be the family of these path of length $n$, ending at level $i$, for $i=0,1,2,3$. 

It is straightforward to prove that $|\mathcal{A}_{2n,0}|=F_{2n-1}$, $|\mathcal{A}_{2n,2}|=F_{2n}$, $|\mathcal{A}_{2n+1,1}|=F_{2n+1}$, and $|\mathcal{A}_{2n+1,3}|=F_{2n}$, with Fibonacci numbers $F_k$.

We also consider similar lattice paths, this time with $-2\le s_i\le 1$ and starting again at the origin. Let 
Let $\mathcal{B}_{n,i}$ be the family of these path of length $n$, ending at level $i$, for $i=-2,-1,0,1$. 
It is again straightforward to prove that $|\mathcal{B}_{2n,0}|=F_{2n+1}$, $|\mathcal{B}_{2n,-2}|=F_{2n}$, $|\mathcal{B}_{2n+1,1}|=F_{2n+1}$, and $|\mathcal{B}_{2n+1,-1}|=F_{2n+2}$.

Thus we have that 
\begin{equation*}
\Big|\bigcup_{i=0}^3\mathcal{A}_{n,i}\Big|=|\mathcal{B}_{n,0}|+|\mathcal{B}_{n,-1}|.
\end{equation*}

Cigler~\cite{cigler} asked for a bijection to explain this fact. An answer was given by Prellberg~\cite{prellberg}.
In this note, I want to link these problems to a tree structure named \emph{Elena} (trees), that I introduced some fifteen years ago~\cite{elena}. The bijection presented in this early paper can also be used here to explain the equality. 

If we want that the paths are in correspondence with trees, we require an even number of steps.

\section{Paths in $\mathcal{A}_{2n,0}$ and height restricted plane trees}

The translation of such a path of length $2n$ into a plane tree of height $\le 3$ (counting edges) is direct and sometimes called \emph{glove bijection}. The following example will be sufficient.

\begin{figure}[h]
	\begin{tikzpicture}[scale=0.5]
	\draw[step=1.cm,black] (0,0) grid (20,3);	
	\node at (-0.6,3) {$3$};	\node at (-0.6,2) {$2$};	\node at (-0.6,0) {$0$};\node at (-0.6,1) {$1$};
	\draw (0,0)to(3,3)to(4,2) to (5,3)to (8,0);
		\draw (8,0)to(9,1)to(10,0) to (11,1)to (12,0)to(14,2)to(15,1)to(16,2)to(17,1)to(18,2)to(20,0);

	\node at (0,0) {$\bullet$};\node at (20,0) {$\bullet$};
	
	\node at (26,4-1) {$\bullet$};
	\node at (23,3-1) {$\bullet$};
	\node at (25,3-1) {$\bullet$};
	\node at (27,3-1) {$\bullet$};
	\node at (29,3-1) {$\bullet$};
	\node at (23,2-1) {$\bullet$};
	\node at (29,2-1) {$\bullet$};
	\node at (28,2-1) {$\bullet$};
	\node at (30,2-1) {$\bullet$};
	\node at (22,1-1) {$\bullet$};
	\node at (24,1-1) {$\bullet$};

	\draw (22,0)to(23,1)to(23,2)to (26,3);\draw (24,0)to(23,1);\draw (25,2)to(26,3)to(27,2);
	\draw (26,3)to(29,2);
	\draw (28,1)to(29,2);
	\draw (29,1)to(29,2);
	\draw (30,1)to(29,2);
	
	\end{tikzpicture}
	\caption{A path of length 20, and the corresponding height restricted plane tree with 11 nodes}
	\end	{figure}
	
	\section{Paths in $\mathcal{B}_{2n,0}$ and Elena trees}
	
	Elenas were introduced in \cite{elena}; they consist of some nodes labelled $\mathbf{a}$, and a sequence of paths
	of various lengths (possibly empty)  emanating from all of them, except for the last one. An example describes this readily:
	\begin{figure}[h]
		\begin{tikzpicture}[scale=0.5]
		
		\draw[ultra thick](0,0)to(4,-4);
		
		\node at (0,0) {$\bullet$};\node at (1,-1) {$\bullet$};\node at (2,-2) {$\bullet$};\node at (3,-3) {$\bullet$};
		\node at (4,-4) {$\bullet$};
		
		\draw(0,0)to(-4.5,-3);\node at (-1.5,-1) {$\bullet$};\node at (-3,-2) {$\bullet$};\node at (-4.5,-3) {$\bullet$};
		\draw(1,-1)to(-2,-2);\draw(1,-1)to(-1.0,-2);\draw(1,-1)to(-3,-5);
		\node at (-2,-2) {$\bullet$};\node at (-1,-2) {$\bullet$};\node at (0,-2) {$\bullet$};\node at (-1,-3) {$\bullet$};
		\node at (-2,-4) {$\bullet$};\node at (-3,-5) {$\bullet$};
		
		\draw(3,-3)to(1,-5);\node at (1,-5) {$\bullet$};\node at (2,-4) {$\bullet$};
		
		\node at (4.5,-3.5){$\mathbf{a}$};
				\node at (3.5,-2.5){$\mathbf{a}$};
						\node at (2.5,-1.5){$\mathbf{a}$};
								\node at (1.5,-0.5){$\mathbf{a}$};
										\node at (0.5,0.5){$\mathbf{a}$};
		\end{tikzpicture}
		\caption{An Elena described by  $\mathbf{a}\mathbf{p}_3\mathbf{a}\mathbf{p}_1\mathbf{p}_1\mathbf{p}_4\mathbf{a}\mathbf{a}\mathbf{p}_2	\mathbf{a}$}
		\end	{figure}
		
		Typically, an Elena can be described by $\mathbf{a}\mathbf{p}_{i_1}\mathbf{p}_{i_2}\dots \mathbf{a}\mathbf{p}_{j_1}\mathbf{p}_{j_2}\dots \mathbf{a} \dots \mathbf{a}$.
		For the set (language) of Elenas, we might write a symbolic expression $\big(\mathbf{a}\mathbf{p}^*\big)^*\mathbf{a}$.
		
		It is perhaps surprising that the paths in $\mathcal{B}_{2n,0}$ are suitable to describe Elenas:
		For each sequence of steps $(2i,0)\to(2i+1,1)\to(2i+2,0)$, we write a symbol $\mathbf{a}$. In Figure~3
		such pairs of steps are depicted in boldface.
		
		Thus, a path can be decomposed as $\mathbf{w}_0\mathbf{a}\mathbf{w}_1\mathbf{a}\dots \mathbf{a}\mathbf{w}_s$, where each $\mathbf{w}$ is a walk from level 0 to level $0$ that 
	``lives'' on levels $0,-1,-2$. Now we add a symbol $\mathbf{a}$ both, to the left and to the right.
	
	What is still left is how such a $\mathbf{w}$ can be interpreted as a sequence of paths: Each return to the level 0 marks the end of a path, and the translation of the sojourns is as follows: 
	
	$\tikz [scale=0.2]\draw[thick] (0,0)to (1,-1)to (2,0);$ corresponds to $\mathbf{p}_1$,
	$\tikz [scale=0.2]\draw[thick] (0,0)to (2,-2)to (4,0);$ corresponds to $\mathbf{p}_2$,
	$\tikz [scale=0.2]\draw[thick] (0,0)to (2,-2)to (3,-1) to (4,-2)to (6,0);$ corresponds to $\mathbf{p}_3$,
		$\tikz [scale=0.2]\draw[thick] (0,0)to (2,-2)to (3,-1) to (4,-2)to (5,-1) to (6,-2)to (8,0);$ corresponds to $\mathbf{p}_4$, and so forth.

		Note that in this way a path of length $2n$ is (bijectively) mapped to an Elena of size (= number of nodes) $n+2$; the Elena consisting only of one node will not be considered. 
		
\begin{figure}[h]
	\begin{tikzpicture}[scale=0.5]
	\draw[step=1.cm,black] (0,1) grid (28,-2);	
		\node at (-0.6,-1) {$-1$};	\node at (-0.6,-2) {$-2$};	\node at (-0.6,0) {$0$};\node at (-0.6,1) {$1$};
			\draw (0,0)to(2,-2)to(3,-1) to (4,-2)to (6,0);
			\draw[ultra thick](6,0)to(7,1)to(8,0);
			\draw (8,0)to(9,-1)to(10,0) to (11,-1)to (12,0)to(14,-2)to(15,-1)to(16,-2)to(17,-1)to(18,-2)to(20,0);
			\draw[ultra thick](20,0)to(21,1)to(22,0)to(23,1)to(24,0);
			\draw (24,0)to(26,-2)to(28,0);
			
			\node at (0,0) {$\bullet$};\node at (28,0) {$\bullet$};
		\end{tikzpicture}
			\caption{A path of length 28, described by $\mathbf{p}_3\mathbf{a}\mathbf{p}_1\mathbf{p}_1\mathbf{p}_4\mathbf{a}\mathbf{a}\mathbf{p}_2$}
			\end{figure}
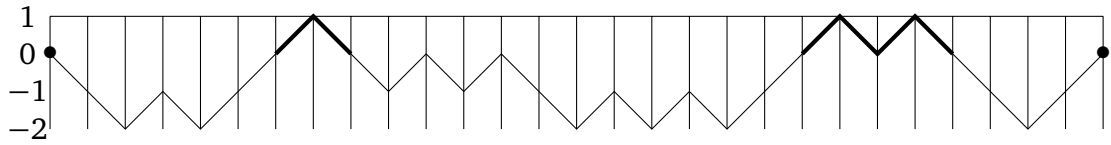

\begin{figure}[h]
	\begin{tikzpicture}[scale=0.5]
	
	\draw[ultra thick](0,0)to(4,-4);
	
	\node at (0,0) {$\bullet$};\node at (1,-1) {$\bullet$};\node at (2,-2) {$\bullet$};\node at (3,-3) {$\bullet$};
	\node at (4,-4) {$\bullet$};
	
	\draw(0,0)to(-4.5,-3);\node at (-1.5,-1) {$\bullet$};\node at (-3,-2) {$\bullet$};\node at (-4.5,-3) {$\bullet$};
	\draw(1,-1)to(-2,-2);\draw(1,-1)to(-1.0,-2);\draw(1,-1)to(-3,-5);
	\node at (-2,-2) {$\bullet$};\node at (-1,-2) {$\bullet$};\node at (0,-2) {$\bullet$};\node at (-1,-3) {$\bullet$};
	\node at (-2,-4) {$\bullet$};\node at (-3,-5) {$\bullet$};
	
	\draw(3,-3)to(1,-5);\node at (1,-5) {$\bullet$};\node at (2,-4) {$\bullet$};
	\end{tikzpicture}
	\caption{The Elena with 16 nodes corresponding to  $(\mathbf{a})\mathbf{p}_3\mathbf{a}\mathbf{p}_1\mathbf{p}_1\mathbf{p}_4\mathbf{a}\mathbf{a}\mathbf{p}_2(\mathbf{a})$}
	\end	{figure}
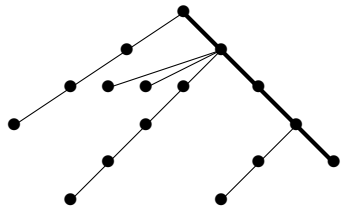

\section{Elenas and height restricted plane trees} 

We will establish a bijection between $\mathcal{B}_{2n,0}$ and 
$\mathcal{A}_{2n,0}\cup\mathcal{A}_{2n,2}$; note, however that the latter set may be replaced by
$\mathcal{A}_{2n+2,0}$, by distinguishing the two cases of the last two steps.

So we would be done once we would know how to map (bijectively)
an Elena of size $n+2$ to a height restricted plane tree of the same size.

 This was documented already in \cite{elena}, but will be repeated here to make this note self contained. The set of operations will be described a sequence of pictures, which require no additional explanation.
 
 We start with our running example of an Elena of size 16 and gradually transform it:
 \begin{figure}[h]
 	\begin{tikzpicture}[scale=0.5]
 	
 	\draw[ultra thick](0,0)to(4,-4);
 	
 	\node at (0,0) {$\bullet$};\node at (1,-1) {$\bullet$};\node at (2,-2) {$\bullet$};\node at (3,-3) {$\bullet$};
 	\node at (4,-4) {$\bullet$};
 	
 	\draw(0,0)to(-4.5,-3);\node at (-1.5,-1) {$\bullet$};\node at (-3,-2) {$\bullet$};\node at (-4.5,-3) {$\bullet$};
 	\draw(1,-1)to(-2,-2);\draw(1,-1)to(-1.0,-2);\draw(1,-1)to(-3,-5);
 	\node at (-2,-2) {$\bullet$};\node at (-1,-2) {$\bullet$};\node at (0,-2) {$\bullet$};\node at (-1,-3) {$\bullet$};
 	\node at (-2,-4) {$\bullet$};\node at (-3,-5) {$\bullet$};
 	
 	\draw(3,-3)to(1,-5);\node at (1,-5) {$\bullet$};\node at (2,-4) {$\bullet$};
 
 \begin{scope}[shift={(10,0)}]
 	\draw[ultra thick](0,0)to(8,0);
 	
 	\node at (0,0) {$\bullet$};\node at (2,0) {$\bullet$};\node at (4,0) {$\bullet$};\node at (6,0) {$\bullet$};
 	\node at (8,0) {$\bullet$};
 	
 	\draw(0,0)to(-2,-3);\node at (-2/3,-1) {$\bullet$};\node at (-4/3,-2) {$\bullet$};\node at (-2,-3) {$\bullet$};
 	\draw(2,0)to(1,-1);\draw(2,0)to(0,-1);\draw(2,0)to(3,-4);
 	\node at (3,-4) {$\bullet$};
 	\node at (2.5,-2) {$\bullet$};\node at (2.25,-1) {$\bullet$};\node at (2.75,-3) {$\bullet$};
 	\node at (0,-1) {$\bullet$};\node at (1,-1) {$\bullet$};
 	
 	\draw(6,0)to(6,-2);\node at (6,-2) {$\bullet$};\node at (6,-1) {$\bullet$};

 \end{scope}
 \begin{scope}[shift={(-4,-8)}]
 \draw[ultra thick](0,0)to(4,1);
 \draw[ultra thick](2,0)to(4,1);
  \draw[ultra thick](4,0)to(4,1);
   \draw[ultra thick](6,0)to(4,1);

 \node at (0,0) {$\bullet$};\node at (2,0) {$\bullet$};\node at (4,0) {$\bullet$};\node at (6,0) {$\bullet$};
 \node at (4,1) {$\bullet$};
 
 \draw(0,0)to(-2,-3);\node at (-2/3,-1) {$\bullet$};\node at (-4/3,-2) {$\bullet$};\node at (-2,-3) {$\bullet$};
 \draw(2,0)to(1,-1);\draw(2,0)to(0,-1);\draw(2,0)to(3,-4);
 \node at (3,-4) {$\bullet$};
 \node at (2.5,-2) {$\bullet$};\node at (2.25,-1) {$\bullet$};\node at (2.75,-3) {$\bullet$};
 \node at (0,-1) {$\bullet$};\node at (1,-1) {$\bullet$};
 
 \draw(6,0)to(6,-2);\node at (6,-2) {$\bullet$};\node at (6,-1) {$\bullet$};
 
 \end{scope}
 \begin{scope}[shift={(6,-8)}]
 \draw[ultra thick](0,0)to(4,1);
 \draw[ultra thick](2,0)to(4,1);
 \draw[ultra thick](4,0)to(4,1);
 \draw[ultra thick](6,0)to(4,1);

 \node at (0,0) {$\bullet$};\node at (2,0) {$\bullet$};\node at (4,0) {$\bullet$};\node at (6,0) {$\bullet$};
 \node at (4,1) {$\bullet$};
 
 \draw(0,0)to(-4/3,-2)to (0,-2);\node at (-2/3,-1) {$\bullet$};\node at (-4/3,-2) {$\bullet$};
 
 \node at (0,-2) {$\bullet$};
 \draw(2,0)to(1,-1);\draw(2,0)to(0,-1);\draw(2,0)to(2,-2)to(4,-2);
 \node at (4,-2) {$\bullet$};
 \node at (2,-2) {$\bullet$};\node at (2,-1) {$\bullet$};\node at (3,-2) {$\bullet$};
 \node at (0,-1) {$\bullet$};\node at (1,-1) {$\bullet$};
 
 \draw(6,0)to(6,-2);\node at (6,-2) {$\bullet$};\node at (6,-1) {$\bullet$};
 
 \end{scope}
 \begin{scope}[shift={(15,-8)}]
 \draw[ultra thick](0,0)to(4,1);
 \draw[ultra thick](2,0)to(4,1);
 \draw[ultra thick](4,0)to(4,1);
 \draw[ultra thick](6,0)to(4,1);

 \node at (0,0) {$\bullet$};\node at (2,0) {$\bullet$};\node at (4,0) {$\bullet$};\node at (6,0) {$\bullet$};
 \node at (4,1) {$\bullet$};
 
 \draw(0,0)to(-4/3,-2);\draw(-2/3,-1)to(0,-2);
 
 \node at (-2/3,-1) {$\bullet$};\node at (-4/3,-2) {$\bullet$};
 
 \node at (0,-2) {$\bullet$};
 \draw(2,0)to(1,-1);\draw(2,0)to(0,-1);\draw(2,0)to(2,-2);
 \draw(2,-1)to(3,-2);
\draw(2,-1)to(4,-2); 
 
 \node at (4,-2) {$\bullet$};
 \node at (2,-2) {$\bullet$};\node at (2,-1) {$\bullet$};\node at (3,-2) {$\bullet$};
 \node at (0,-1) {$\bullet$};\node at (1,-1) {$\bullet$};
 
 \draw(6,0)to(6,-2);\node at (6,-2) {$\bullet$};\node at (6,-1) {$\bullet$};
 
 \end{scope}
 
 	\end{tikzpicture}
 
 \caption{Transforming an Elena into a height restricted plane tree}
 
 		\end{figure}
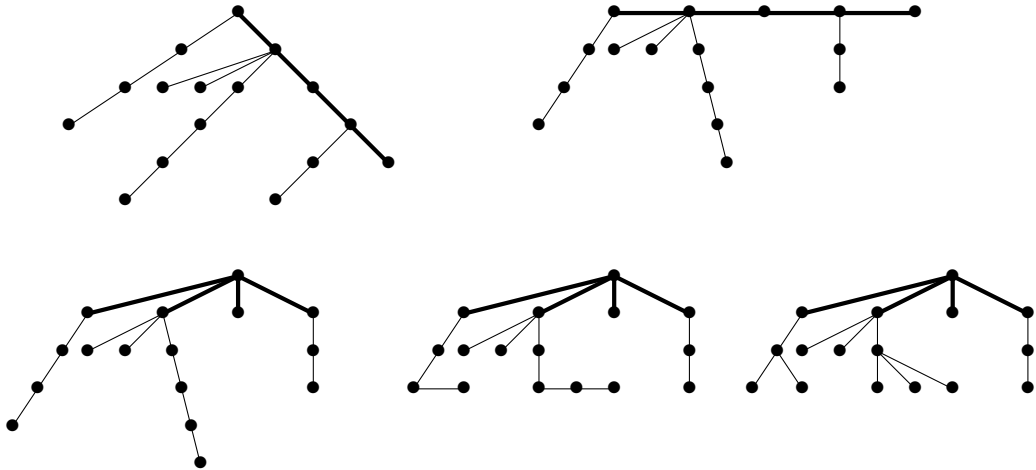

 \section{Paths with an odd number of steps}
 	
 	Let us consider $\mathcal{B}_{2n-1,-1}$, enumerated by $F_{2n}$. If we augment one up-step at the end,
 	we have Elenas, but with the special property that the last group of paths is non-empty.
 	
 	One the other hand, if we consider $\mathcal{A}_{2n-1,1}\cup\mathcal{A}_{2n-1,3}$, which is equivalent to
 	$\mathcal{A}_{2n,2}$, then we augment it with 2 down-steps. The resulting height restricted tree has the property that the rightmost leaf is on a level $\ge2$. 
 	
 	A short reflection convinces us that the bijection described earlier also works bijectively on the two respective subclasses.

\end{document}